# Trees, permutations and the tangent function

Ross Street[1]

27 July 2001

## §0. Introduction

The Macquarie University Mathematics Department offers four-week Summer Vacation Scholarships to bright undergraduate students. The two Scholars I supervised at the beginning of this year were Ryan Crompton and Tam Pham — this is a talk on some of their research.

Towards the end of last year, my postdoctoral research associate William Joyce gave a seminar on an application of Category Theory to Physics [Jo]. He made use of a particular kind of mathematical tree and produced a formula for iteratively calculating the number of these trees with a given number of leaves. As we allow the number of leaves to increase we obtain a sequence of numbers.

Now, there is a Web Page [Sl1] operated by N.J.A. Sloane which can tell you, from typing in the first few terms of a sequence, whether that sequence has occurred somewhere else in Mathematics. My postgraduate student Daniel Steffen traced this down and found, to our huge surprise, that the sequence was related to the tangent function $\tan x$. Ryan and Tam searched out what was known about this connection and discovered some apparently new results.

We all found this a lot of fun and I hope you will too.

## §1. Permutations and combinations

We begin with some familiar combinatorial ideas. A *permutation* of a list of distinct numbers is a list of the same numbers in a different order; this second list must involve all the numbers of the original list without repetition. For example, one permutation of the list 2 3 6 7 9 is 6 3 9 2 7. The total number of permutations of any list

$$a_1\ a_2\ \ldots\ a_n$$

of length n is

$$n! = n \times (n-1) \times (n-2) \times \ldots \times 2 \times 1\ ;$$

the exclamation mark is called *factorial* in this context.

For our purposes, the original list $a_1\ a_2\ \ldots\ a_n$ will be in increasing order:

$$a_1 < a_2 < \ \ldots\ < a_n\ .$$

A choice of some of the numbers $a_1\ a_2\ \ldots\ a_n$ can always be written in increasing order too: such a choice is called a *combination* from the given list. The number of combinations of length r from an original list of length n is

[1] A talk for the HSC Mathematics Talented Students Day at Macquarie University on 27 July 2001.

$$^{n}C_r = \frac{n!}{r!\,(n-r)!}\,;$$

these numbers $^{n}C_r$ are called *binomial coefficients.*

The binomial coefficients can be constructed by a triangular process noticed by Blaise Pascal (1623 – 1662). Each row has a 1 on each end and the other entries are the sum of the two above it. The r-th entry in row n of *Figure 1* is $^{n}C_r$ (we start counting at $n = 0$ and $r = 0$).

| | | | | | | | | | | | | | | | |
|---|---|---|---|---|---|---|---|---|---|---|---|---|---|---|---|
| n=0 | | | | | | | | 1 | | | | | | | |
| 1 | | | | | | | 1 | | 1 | | | | | | |
| 2 | | | | | | 1 | | 2 | | 1 | | | | | |
| 3 | | | | | 1 | | 3 | | 3 | | 1 | | | | |
| 4 | | | | 1 | | 4 | | 6 | | 4 | | 1 | | | |
| 5 | | | 1 | | 5 | | 10 | | 10 | | 5 | | 1 | | |
| 6 | | 1 | | 6 | | 15 | | 20 | | 15 | | 6 | | 1 | |
| 7 | 1 | | 7 | | 21 | | 35 | | 35 | | 21 | | 7 | | 1 |

*Figure 1*

## §2. Polynomial approximations

If we know that a function $f(x)$ is well behaved and we know the values

$$f(0),\ f'(0),\ f''(0),\ \ldots$$

of the function and its repeated derivatives at 0, we know a lot about the function. For example, suppose that we suspect our function is a cubic

$$f(x) = a + bx + cx^2 + dx^3,$$

then $f'(x) = b + 2cx + 3dx^2$, $f''(x) = 2c + 6dx$, $f'''(x) = 6d$ ; so we find

$$a = f(0)\,,\ b = f'(0)\,,\ c = \tfrac{1}{2!}f''(0)\,,\ d = \tfrac{1}{3!}f'''(0)\ .$$

You will be familiar with the fact that $\sin x$ is well approximated by the polynomial $x$ for small $x$. For $f(x) = \sin x$ we can easily calculate that $f(0) = \sin 0 = 0$, $f'(0) = \cos 0 = 1$, $f''(0) = -\sin 0 = 0$, $f'''(0) = -\cos 0 = -1$. Of course, we do not really expect $\sin x$ to be a cubic, but if we use the formulas above for a, b, c, d, we obtain the cubic $x - \frac{1}{3!}x^3$. Indeed, you can check on your calculator that this gives a better approximation to $\sin x$ than $x$ for small $x$. Taking this further, we obtain a polynomial

$$x - \tfrac{1}{3!}x^3 + \tfrac{1}{5!}x^5 - \tfrac{1}{7!}x^7 + \tfrac{1}{9!}x^9 - \ldots + (-1)^n \tfrac{1}{(2n+1)!}x^{2n+1}$$

which can be used to approximate $\sin x$ for all $x$ by taking $n$ large enough. Polynomials obtained from functions in this way, by using the values of the derivatives at 0, are called *Maclaurin* or *Taylor* polynomials.

In the case of $y = \sin x$, finding the derivatives at 0 is made easy by the fact that the derivatives start repeating after four steps: $y'''' = \sin x$. Even before that we see that

$$y'' = -y .$$

Whenever we have an expression for a higher derivative of a function in terms of its earlier derivatives (what is called a *differential equation*), we can use it to find the Taylor polynomials.

Let us now focus our attention on the tangent function $y = \tan x$. We can record that $y(0) = 0$. We have

$$y' = \sec^2 x = 1 + \tan^2 x = 1 + y^2$$

yielding $y'(0) = 1$ and the differential equation $y' = 1 + y^2$. Differentiating again we obtain another differential equation

$$y'' = 2\,y\,y',$$

yielding $y''(0) = 0$. Continuing, we see that

$$\begin{aligned} y''' &= 2\,y'y' + 2\,y\,y'', \\ y^{(4)} &= 6\,y'y'' + 2\,y\,y''', \\ y^{(5)} &= 6\,y''y'' + 8\,y'\,y''' + 2\,y\,y^{(4)}, \\ y^{(6)} &= 20\,y''y''' + 10\,y'\,y^{(4)} + 2\,y\,y^{(5)}, \\ y^{(7)} &= 20\,y'''\,y''' + 30\,y''\,y^{(4)} + 12\,y'\,y^{(5)} + 2\,y\,y^{(6)}, \end{aligned}$$

yielding $y'''(0) = 2$, $y^{(4)}(0) = 0$, $y^{(5)}(0) = 16$, $y^{(6)}(0) = 0$, $y^{(7)}(0) = 272$. It can be checked that the Taylor polynomial

$$x + 2\,\frac{x^3}{3!} + 16\,\frac{x^5}{5!} + 272\,\frac{x^7}{7!}$$

is quite a good approximation to $\tan x$; you can obtain a better approximation by finding higher $y^{(n)}(0)$. Notice that it is fairly predictible that $y^{(n)}(0) = 0$ for all even $n$.

## §3. Joyce trees

Graphs of various kinds are used throughout Mathematics and its applications. For example, graphs can be used to organize both computer hardware design and computer programs. Graphs have nodes (represented by points) and edges (represented by curves – straight lines if possible) connecting various nodes.

A *tree* is a graph which is connected (you can get from one node to any other by a path of edges) and loop free (you cannot get back to the same node once you set out on a path of edges, without backtracking). We are only interested in trees (see Figure 2) for which a particular node is selected and called the *root*. We then draw trees on a plane piece of paper with the root at the bottom, thought of as at level 0. We wish to record the fact that certain nodes are at the same level. For example, in Figure 2, there are two nodes at level 1, three at level 2, and four at level 3. Nodes with no edges connected above them are called *leaves*. The tree in Figure 2 has six leaves.

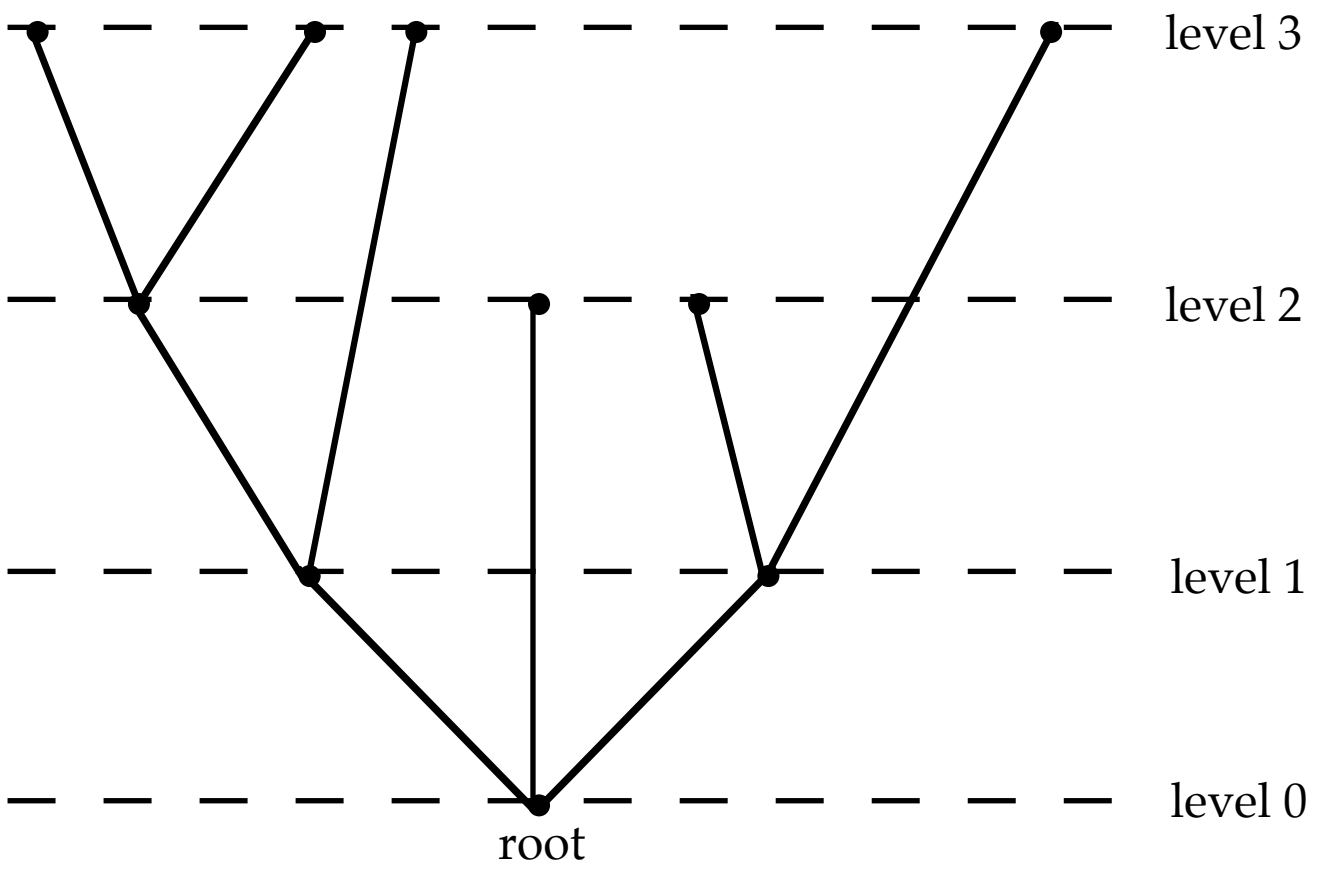

*Figure 2*

A tree is called *binary* when each node is either a leaf or has precisely two edges above connected to it. Figure 2 is not a binary tree since the root has three edges connected above it. A *Joyce tree* is a binary tree for which no two nodes have the same level and all levels, up to that of the top leaf, have a node. An example is provided in Figure 3.

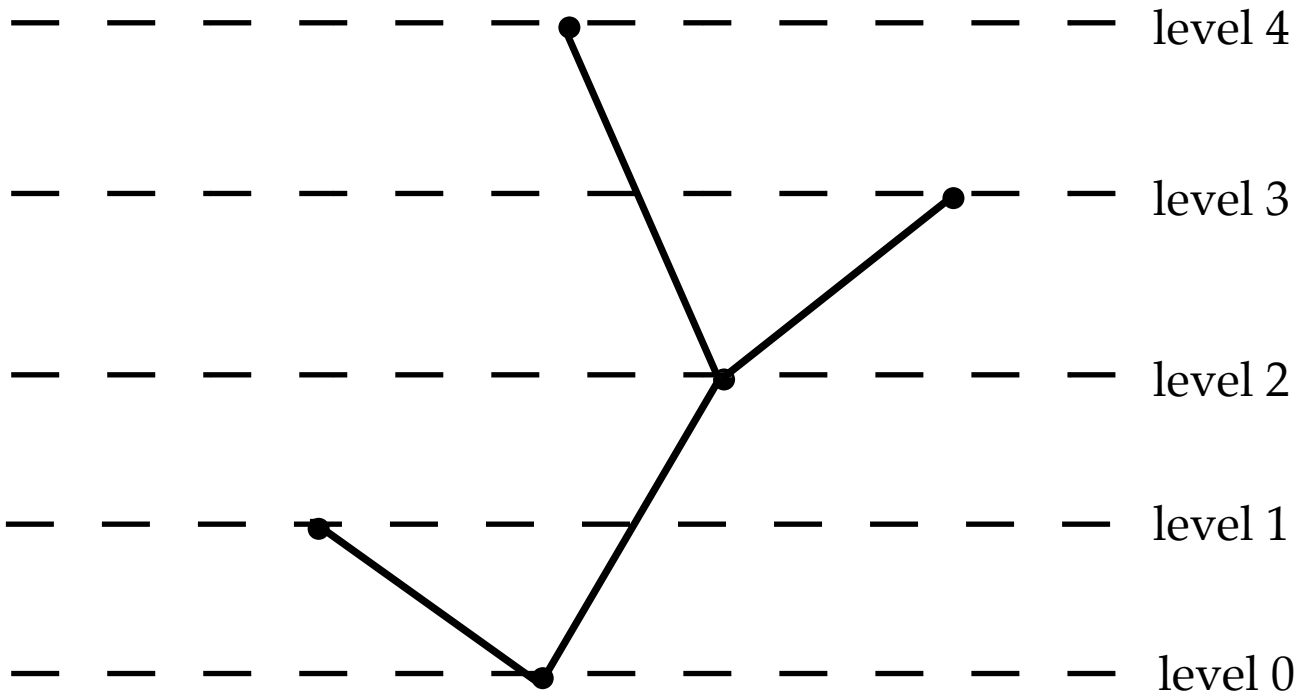

*Figure 3*

A fact we need to observe about a Joyce tree is that there is a simple relationship between the number of leaves and the number of nodes: if there are $m$ leaves then there are $2m-1$ nodes altogether. In Figure 3, there are $m=3$ leaves and $2 . 3-1=5$ nodes. To prove this in general notice that, if we have a Joyce tree with more than one leaf, look at the highest level where there is a node which is not a leaf. The two edges above and connected to this node must connect the node to two leaves. Removal of those two edges and the two leaves creates a new leaf. So we have reduced the number of leaves by one and the number of nodes by two. What remains is still a Joyce tree. The process continues until we have a single node. The fact we want now follows by induction.

Let $J_n$ be the number of Joyce trees with $n$ nodes. One thing we have seen is that $J_n = 0$ for $n$ even.

## §4. Tremolo permutations

A permutation of an increasing list of numbers is called *tremolo* when consecutive differences between consecutive numbers in the permuted list have opposite sign. For example, the permutation 6 3 5 2 4 1 of 1 2 3 4 5 6 is tremolo since 6 – 3 is positive, 3 – 5 is negative, 5 – 2 is positive, 2 – 4 is negative, 4 – 1 is positive. It reminds me of playing tremolo on a mandolin: down-up-down-up-down.

There is a one-to-one correspondence between Joyce trees with n nodes and tremolo permutations of 0 1 2 3 . . . n n+1 which begin with 1 and end with 0. To see this, label the nodes of the tree with numbers obtained by adding 2 to the level of the node. Now we "read" the tree from left to right and from top to bottom taking note of the numbers labelling the nodes; this gives a tremolo permutation of 2 3 4 . . . n n+1. Put the 1 at the front and 0 at the end. For example, for the tree in Figure 3 we obtain the labelling of nodes as in Figure 4. Reading the tree gives the list 3 2 6 4 5. So the desired tremolo permutation of 0 1 2 3 4 5 6 is 1 3 2 6 4 5 0. Notice that the numbers in even positions (the second, fourth, sixth, . . ) in the permutation always correspond to leaves. Notice too that our insistence on a 1 at the start and 0 at the end means there can be no tremolo permutations of 0 1 2 3 . . . n n+1 for n even.

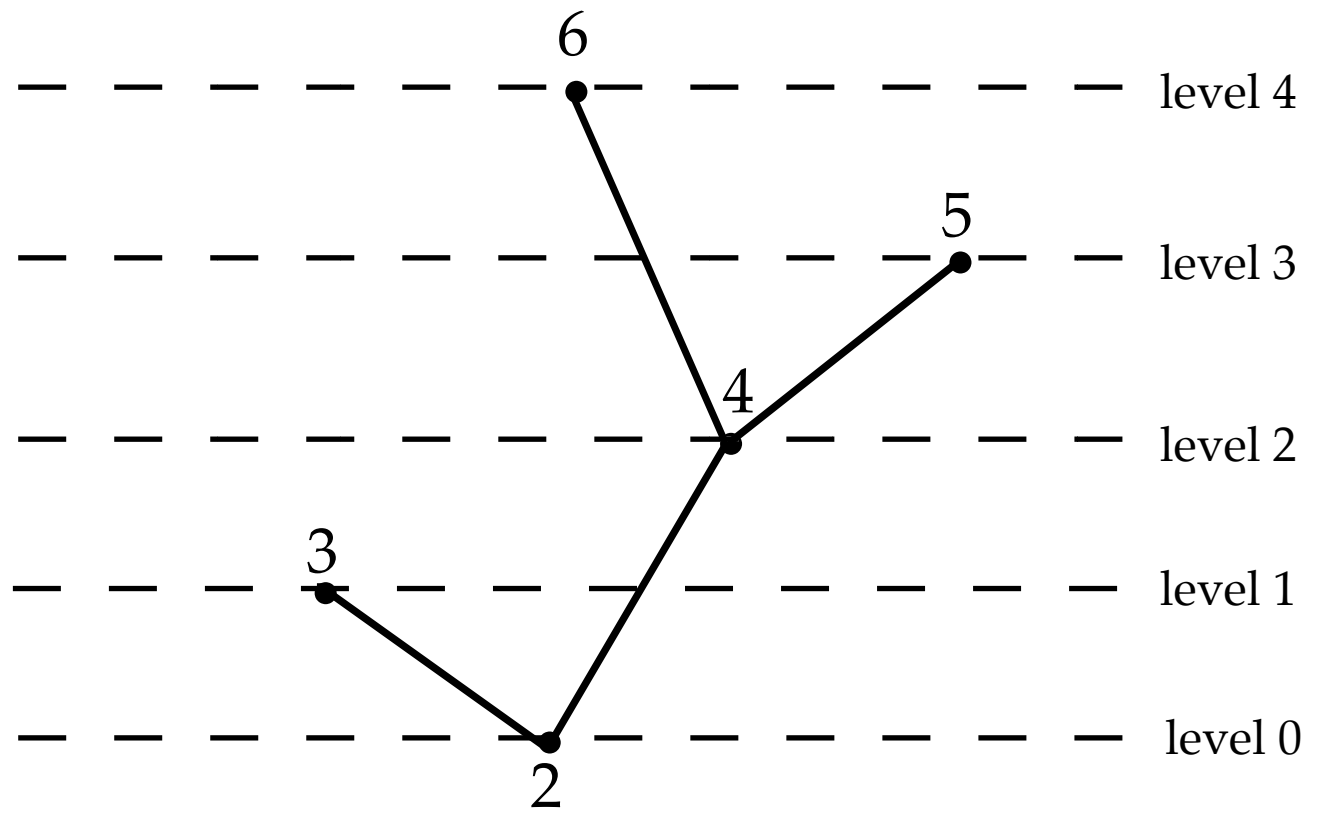

*Figure 4*

Clearly what this shows is that $J_n$ is equal to the number of tremolo permutations of any list of n+2 increasing numbers which have the smallest number at the end and the second smallest at the beginning. Equally, $J_n$ is the number of tremolo permutations of any list of n+2 increasing numbers which have the smallest number at the beginning and second smallest at the end (just read the permutation backwards).

We shall now find a formula for recursively determining the numbers $J_n$ ; we shall express $J_{n+1}$ in terms of the $J_m$ with $m < n$. Consider any tremolo permutation of 0 1 2 3 . . . n+2 beginning with 1 and ending with 0:

$$1\overbrace{****************}^{n+1}0.$$

How is such a permutation made up? First we must choose one of the n starred positions to place the number 2; let us say we choose the (m+1)-th star.

$$1\overbrace{*****}^{m}\,2\,\overbrace{**********}^{n-m}\,0.$$

Now we need to choose m numbers from the list 3 4 5 . . . n+2 of n remaining numbers; this choice can be made in ${}^nC_m$ ways. Using the chosen numbers to insert between 1 and 2, we need to choose a tremolo permutation starting at 1 and ending at 2; this choice can be made in $J_m$ ways. Using the remaining n–m numbers, we need to choose a tremolo permutation starting at 2 and ending at 0; this can be done in $J_{n-m}$ ways. This means, once we fix the position of 2 at the (m+1)-th star, there are ${}^nC_m\, J_m\, J_{n-m}$ possibilities. It follows that

$$J_{n+1} = \sum_{m=0}^{n} {}^nC_m\, J_m\, J_{n-m} .$$

This formula can be used to calculate the $J_n$ given that we know $J_0 = 0$ and $J_1 = 1$. For then,

$$J_2 = {}^1C_0\, J_0\, J_1 + {}^1C_1\, J_1\, J_0 = 0\ , \quad J_3 = {}^2C_0\, J_0\, J_2 + {}^2C_1\, J_1\, J_1 + {}^2C_2\, J_2\, J_0 = 0 + 2 + 0 = 2\,,$$

by which time it becomes clear why $J_n = 0$ for n even and why we need only worry about m odd in the summation, so

$$J_5 = {}^4C_1\, J_1\, J_3 + {}^4C_3\, J_3\, J_1 = 4\times 2 + 4\times 2 = 16,$$

$$J_7 = {}^6C_1\, J_1\, J_5 + {}^6C_3\, J_3\, J_3 + {}^6C_5\, J_5\, J_1 = 6\times 1\times 16 + 20\times 2\times 2 + 6\times 16\times 1 = 272,$$

$J_9 = 7936$, $J_{11} = 353792$, $J_{13} = 22368256$, $J_{15} = 1903757312$, $J_{17} = 209865342976$,
$J_{19} = 29088885112832$, $J_{21} = 4951498053124096$, $J_{23} = 1015423886506852352$,
$J_{25} = 246921480190207983616$, $J_{27} = 70251601603943959887872$,
$J_{29} = 23119184187809597841473536$, . . . .

Comparing with Section 2, we have empirical evidence that the $J_n$ are the coefficients in the Taylor polynomials of $\tan x$. In Section 5 we shall show how to prove it completely.

## §5. Back to tan x

We know that $y = \tan x$ satisfies the differential equation $y' = 1 + y^2$ . Let us see whether we can use this to find the numbers $t_n$ in order for the expression

$$y = \frac{t_0}{0!}x^0 + \frac{t_1}{1!}x^1 + \frac{t_2}{2!}x^2 + \frac{t_3}{3!}x^3 + \frac{t_4}{4!}x^4 + \frac{t_5}{5!}x^5 + \frac{t_6}{6!}x^6 + \ldots$$

to satisfy the condition $y' = 1 + y^2$. Differentiating the sum for y, we obtain

$$y' = \frac{t_1}{0!}x^0 + \frac{t_2}{1!}x^1 + \frac{t_3}{2!}x^2 + \frac{t_4}{3!}x^3 + \frac{t_5}{4!}x^4 + \frac{t_6}{5!}x^5 + \frac{t_7}{6!}x^6 + \ldots .$$

Squaring the sum for y, we obtain

$$y^2 = \frac{t_0\,t_0}{0!\,0!}x^0 + \left(\frac{t_0\,t_1}{0!\,1!} + \frac{t_1\,t_0}{1!\,0!}\right)x^1 + \left(\frac{t_0\,t_2}{0!\,2!} + \frac{t_1\,t_1}{1!\,1!} + \frac{t_2\,t_0}{2!\,0!}\right)x^2 + \left(\frac{t_0\,t_3}{0!\,3!} + \frac{t_1\,t_2}{1!\,2!} + \frac{t_2\,t_1}{2!\,1!} + \frac{t_3\,t_0}{3!\,0!}\right)x^3 + \ldots.$$

So, we can equate coefficients of the powers of $x$ in the equation $y' = 1 + y^2$:

$$\frac{t_1}{0!} = 1 + \frac{t_0\,t_0}{0!\,0!}, \quad \frac{t_2}{1!} = \frac{t_0\,t_1}{0!\,1!} + \frac{t_1\,t_0}{1!\,0!}, \quad \frac{t_3}{2!} = \frac{t_0\,t_2}{0!\,2!} + \frac{t_1\,t_1}{1!\,1!} + \frac{t_2\,t_0}{2!\,0!}, \quad \frac{t_4}{3!} = \frac{t_0\,t_3}{0!\,3!} + \frac{t_1\,t_2}{1!\,2!} + \frac{t_2\,t_1}{2!\,1!} + \frac{t_3\,t_0}{3!\,0!},$$

and so on, from which we see that

$$t_{n+1} = \sum_{m=0}^{n} {}^{n}C_m\, t_m\, t_{n-m}.$$

This is precisely the same formula satisfied by the $J_n$. We have already seen that this formula, together with the fact that $t_0 = 0$ and $t_1 = 1$, determines the numbers $t_n$. This proves that $J_n = t_n$ for all $n$.

## §6. Ox ploughing

There is something like a Pascal triangle which can be used to generate the numbers $J_n$. The triangle contains many more numbers than the $J_n$. Let us write ${}^{n}B_m$ for the m-th entry in the n-th row (this time starting with n = 1 and m = 1).

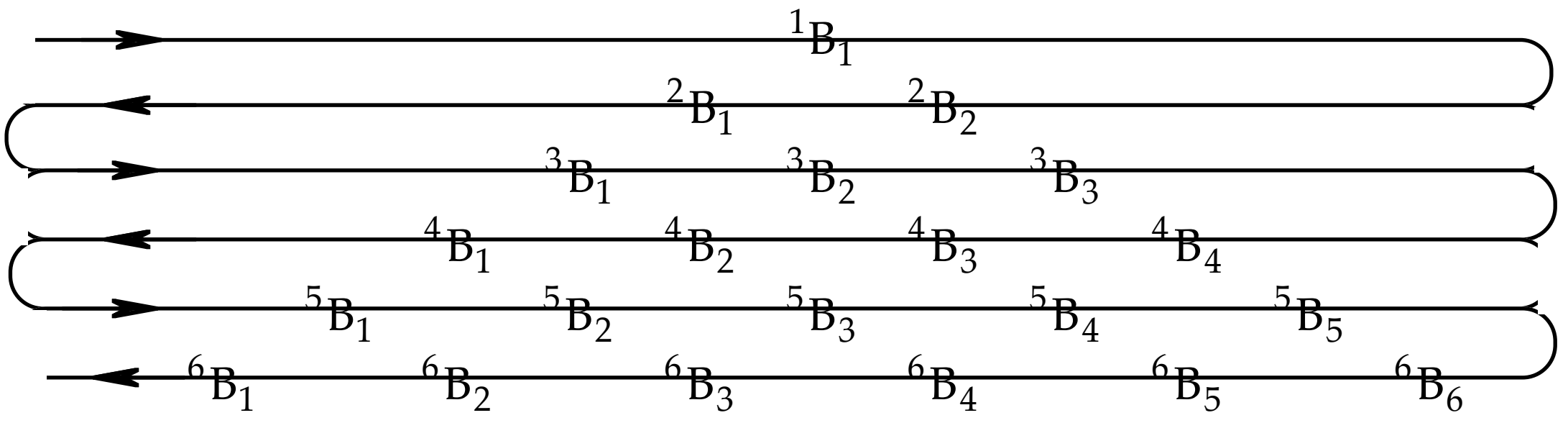

*Figure 5*

In traversing this triangle we move as if we were ploughing rows in a field with an ox; this is called the *Boustrophedon order.* Using this order, we start with ${}^{1}B_1 = 1$, however, we begin each new row with 0 (so that ${}^{2}B_2 = {}^{3}B_1 = {}^{4}B_4 = {}^{5}B_1 = \ldots = 0$ ), and we obtain all other ${}^{n}B_m$ by adding the one before it to the one above and between those two. To be precise, we have

$${}^{1}B_1 = 1, \quad {}^{2k}B_{2k} = {}^{2k+1}B_1 = 0, \quad {}^{2k}B_m = {}^{2k}B_{m+1} + {}^{2k-1}B_m, \quad {}^{2k+1}B_m = {}^{2k+1}B_{m\pm1} + {}^{2k}B_{m-1}.$$

This process leads to the Boustrophedon triangle shown in Figure 6. We can easily see from this construction that each ${}^{n}B_m$ with $n$ even is obtained by adding all the entries to the right of it in the row above, and that each ${}^{n}B_m$ with $n$ odd is obtained by adding all the entries to the left of it in the row above; more formally:

$$ {}^{2k}B_m = \sum_{r=m}^{2k-1} {}^{2k-1}B_r \ , \quad {}^{2k+1}B_m = \sum_{r=1}^{m-1} {}^{2k}B_r \ . $$

The thing to notice about Figure 6 is that the left hand side of the triangle again gives the numbers $J_n$. This empirically verifies the simple identity

$$ {}^{n+1}B_1 = J_n \quad \text{for } n > 0 . $$

This identity can actually be proved rigorously. It follows from a much nicer observation of Ms Tam Pham which identifies all the entries in the Boustrophedon triangle.

```
                              1
                           1     0
                        0     1     1
                     2     2     1     0
                  0     2     4     5     5
              16    16    14    10     5     0
            0    16    32    46    56    61    61
        272   272   256   224   178   122    61     0
      0   272   544   800  1024  1202  1324  1385  1385
 7936  7936  7664  7120  6320  5296  4094  2770  1385     0
```

Boustrophedon triangle

*Figure 6*

**Proposition** (Tam Pham [PC]) *The entry* $^{n}B_m$ *in the Boustrophedon triangle is the number of tremolo permutations of* $0\,1 \ldots n$ *that begin with* $m$ *and end with* $0$.

**Proof** We shall verify that, with this interpretation of $^{n}B_m$ in terms of tremolos, the Boustrophedon construction is satisfied. The idea is fully understandable from looking at the special case $^{8}B_4 = {}^{7}B_4 + {}^{7}B_5 + {}^{7}B_6 + {}^{7}B_7$ . So consider a tremolo permutation

$$ 4 * * * * * * * 0 $$

of 0 1 2 3 4 5 6 7 8. Because 8 is even, in order to tremelo to 0 in 8 steps, the first star $*$ must be greater than 4. So we can leave out the 4 and see that the number of such tremelos is equal to the number of tremelos of 0 1 2 3 5 6 7 8 which start with a number greater than 4 and end with 0. By renumbering the 5 6 7 8 as 4 5 6 7, we see that this is the same as the number of tremelos of 0 1 2 3 4 5 6 7 which start with a number greater than 3 and end with 0. This is the same as the number of tremolos of 0 1 2 3 4 5 6 7 ending in 0 and starting with 4 or 5 or 6 or 7. So the tremolo interpretation of the $^{n}B_m$ does indeed satisfy $^{8}B_4 = {}^{7}B_4 + {}^{7}B_5 + {}^{7}B_6 + {}^{7}B_7$, quod erat demonstrandum (q. e. d.).

**Exercise** Show that the right-hand side of the Boustrophedon Triangle is the sequence of Taylor coefficients for the function $y = \sec x$.

## References and extra reading

Centre of Australian Category Theory
Macquarie University
N. S. W. 2109
AUSTRALIA
email: <street@math.mq.edu.au>
homepage: <http://www.math.mq.edu.au/~street/>